\documentclass[12pt]{amsart}
\usepackage{amsmath,amsthm,amsfonts,amssymb,latexsym}
\usepackage{tikz-cd}
\usepackage{hyperref}
\usepackage{ulem}
\usepackage{pgf,tikz}
\usepackage{mathrsfs}
\usetikzlibrary{arrows}

\usepackage{enumerate}
\usepackage[shortlabels]{enumitem}
\usepackage{color}
\usepackage{xcolor} 
\begin{document}

\theoremstyle{plain}

\newtheorem{thm}{Theorem}[section]
\newtheorem{rem}[thm]{Remark}
\newtheorem{corollary}[thm]{Corollary}
\newtheorem{lemma}[thm]{Lemma}
\newtheorem{theorem}[thm]{Theorem}
\newtheorem{proposition}[thm]{Proposition}
\theoremstyle{definition}
\newtheorem{remark}{Remark}
\newtheorem*{definition}{Definition}
\newtheorem*{notation}{Notation}
\newtheorem{example}[thm]{Example}
\newtheorem*{ThmA}{Theorem A}
\newtheorem*{ThmB}{Theorem B}
\newtheorem*{ThmC}{Theorem C}
\newtheorem*{PropB}{Proposition B}
\newtheorem*{CorB}{Corollary B}
\newtheorem*{CorC}{Corollary C}
\newtheorem*{CorD}{Corollary D}
\newtheorem*{ConjectureB}{Conjecture B}
\newenvironment{enumeratei}{\begin{enumerate}[\upshape (a)]}
    {\end{enumerate}}

\newenvironment{Enumeratei}{\begin{enumerate}[\upshape (A)]}
    {\end{enumerate}}

\def\irr#1{{\rm Irr}(#1)}
\def\cent#1#2{{\bf C}_{#1}(#2)}
\def\pow#1#2{{\mathcal{P}}_{#1}(#2)}
\def\syl#1#2{{\rm Syl}_#1(#2)}
\def\hall#1#2{{\rm Hall}_#1(#2)}
\def\nor{\trianglelefteq\,}
\def\oh#1#2{{\bf O}_{#1}(#2)}
\def\zent#1{{\bf Z}(#1)}
\def\sbs{\subseteq}
\def\gen#1{\langle#1\rangle}
\def\aut#1{{\rm Aut}(#1)}
\def\out#1{{\rm Out}(#1)}
\def\gv#1{{\rm Van}(#1)}
\def\fit#1{{\bf F}(#1)}
\def\frat#1{{\bf \Phi}(#1)}
\def\gammav#1{{\Gamma}(#1)}
\newcommand{\p}{{\mathbb P}}
\newcommand{\N}{{\mathbb N}}
\newcommand{\F}{{\mathbb F}}
\def\fitd#1{{\bf F}_{2}(#1)}
\def\irr#1{{\rm Irr}(#1)}
\def\dl#1{{\rm dl}(#1)}
\def\h#1{{\rm h}(#1)}
\def\ibr#1#2{{\rm IBr}_#1(#2)}
\def\cs#1{{\rm cs}(#1)}
\def\m#1{{\rm m}(#1)}
\def\n#1{{\rm n}(#1)}
\def\cent#1#2{{\bf C}_{#1}(#2)}
\def\hall#1#2{{\rm Hall}_#1(#2)}
\def\syl#1#2{{\rm Syl}_#1(#2)}
\def\nor{\trianglelefteq\,}
\def\norm#1#2{{\bf N}_{#1}(#2)}
\def\oh#1#2{{\bf O}_{#1}(#2)}
\def\Oh#1#2{{\bf O}^{#1}(#2)}
\def\zent#1{{\bf Z}(#1)}
\def\sbs{\subseteq}
\def\gen#1{\langle#1\rangle}
\def\aut#1{{\rm Aut}(#1)}
\def\gal#1{{\rm Gal}(#1)}
\def\alt#1{{\rm Alt}(#1)}
\def\sym#1{{\rm Sym}(#1)}
\def\out#1{{\rm Out}(#1)}
\def\gv#1{{\rm Van}(#1)}
\def\fit#1{{\bf F}(#1)}
\def\lay#1{{\bf E}(#1)}
\def\fitg#1{{\bf F^*}(#1)}

\def\GF#1{{\rm GF}(#1)}
\def\SL#1{{\rm SL}_{2}(#1)}
\def\PSL#1{{\rm PSL}_{2}(#1)}

\def\gammav#1{{\Gamma}(#1)}
\def\V#1{{\rm V}(#1)}
\def\E#1{{\rm E}(#1)}
\def\b#1{\overline{#1}}

 \def\sl#1#2{{\rm SL}_{#1}(#2)}
 \def\gl#1#2{{\rm GL}_{#1}(#2)}
 \def\cl#1#2{{\rm cl}_{#1}(#2)}
\def\Z{{\mathbb{Z}}}
\def\C{{\Bbb C}}
\def\Q{{\Bbb Q}}
\def\inv{^{-1}}
\def\irr#1{{\rm Irr}(#1)}
\def\irrv#1{{\rm Irr}_{\Bbb R}(#1)}
 \def\irrk#1{{\rm Irr}_{ {\rm rv}, K}(#1)}
 \def\irrc#1{{\rm Irr}_{C}(#1)}
  \def\irrf#1{{\rm Irr}_{\mathfrak{F}'}(#1)}
   \def\ext#1{{\rm Ext}(#1)}
   \def\irrh#1{{\rm Irr}_{H}(#1)}
  \def\re#1{{\rm Re}(#1)}
  \def\csrv#1{{\rm cs}_{\rm rv}(#1)}
   \def\clrv#1{{\rm Cl}_{\rm rv}(#1)}
     \def\clk#1{{\rm Cl}_{{\rm rv}, K}(#1)}
  \def\bip#1{{\rm B}_{p'}(#1)}
  \def\irra#1{{\rm Irr}_A(#1)}
   \def\irrs#1{{\rm Irr}_\sigma(#1)}
   \def\irrp#1{{\rm Irr}_{p'}(#1)}
 \def\cdrv#1{{\rm cd}_{\rm rv}(#1)}
 \def\bip#1{{\rm B}_{p'}(#1)}
\def\cdrv#1{{\rm cd}_{\rm rv}(#1)}
\def\cd#1{{\rm cd}(#1)}
\def\irrat#1{{\rm Irr}_{\rm rat}(#1)}
\def\cdrat#1{{\rm cd}_{\rm rat}(#1)}
\def \c#1{{\cal #1}}
\def\cent#1#2{{\bf C}_{#1}(#2)}
\def\syl#1#2{{\rm Syl}_#1(#2)}
\def\oh#1#2{{\bf O}_{#1}(#2)}
\def\Oh#1#2{{\bf O}^{#1}(#2)}
\def\zent#1{{\bf Z}(#1)}
\def\det#1{{\rm det}(#1)}
\def\ker#1{{\rm ker}(#1)}
\def\norm#1#2{{\bf N}_{#1}(#2)}
\def\alt#1{{\rm Alt}(#1)}
\def\iitem#1{\goodbreak\par\noindent{\bf #1}}
    \def \mod#1{\, {\rm mod} \, #1 \, }
\def\sbs{\subseteq}

\def\Char{\rm Char}
\def\Irr{\rm Irr}
\def\Ext{\rm Ext}
\def\Syl{\rm Syl}

\def\Min#1#2{{\rm Min}_{#1}(#2)}
\def\Mcd#1#2{{\rm Mcd}_{#1}(#2)}
\def\Lcd#1#2{{\rm Lcd}_{#1}(#2)}
\def\lin#1#2{{\rm Lin}_{#1}(#2)}
\def\Lin#1{\rm Lin(#1)}
\def\csg#1#2{{\rm cs}_{#1}(#2)}

\setlist[itemize]{font=\color{itemizecolor}}
\colorlet{itemizecolor}{.}
\setlist[enumerate]{font=\color{enumeratecolor}}
\colorlet{enumeratecolor}{.}
\setlist[description]{font=\bfseries\color{descriptioncolor}}
\colorlet{descriptioncolor}{.}

\def\cE{\bar{\rm E}}

\def \nq{\mathfrak{N}_q}

\marginparsep-0.5cm

\renewcommand{\thefootnote}{\fnsymbol{footnote}}
\footnotesep6.5pt

\title[]{Linking conjugacy classes and minimal invariant characters of normal subgroups}
\maketitle

\bigskip

\begin{center}
M.J. Felipe, I. Gilabert, L. Sanus
\end{center}

\begin{abstract} 
Let $G$ be a finite group and $N$ a normal subgroup of $G$. We report on recent results concerning minimal $G$-invariant characters of $N$ (which are the sums of the characters on each orbit of the action of $G$ by conjugation on $\irr N$) and their influence on the structure of $N$, as well as their relationship to the $G$-conjugacy classes of $N$.
\end{abstract}

\bigskip


\bigskip


\bigskip

\section{Introduction}

All groups considered here are finite. A classical yet still compelling observation in the study of finite groups is the parallelism between results concerning conjugacy classes and irreducible characters; in particular, the not completely explained duality between how the sizes of conjugacy classes and the degrees of irreducible characters affect the structure of a group.

Let $G$ be a finite group and let $N$ be a normal subgroup of $G$. During the last few decades, several research papers have moved the focus to the study of \textit{$G$-conjugacy classes of $N$} (recall that they are the classes in $G$ constituted by elements of $N$), and to how they affect the structure of the normal subgroup.

It is known that $G$ acts by conjugation on both the conjugacy classes of $N$ (determining its $G$-conjugacy classes) and on the characters of $N$. A character $\mu$ of $N$ is said to be $G$-invariant if $\mu^g=\mu$ for all $g \in G$. In particular, $\mu$ is said to be \textit{minimal $G$-invariant} if it cannot be expressed as sum of two $G$-invariant characters. Suppose $\theta$ is an irreducible character of $N$ and let $\{\theta^{g_1},\hdots,\theta^{g_{t_\theta}}\}$ be the distinct conjugates of $\theta$ under the action of $G$. Then $\widehat{\theta}=\theta^{g_1}+\cdots+\theta^{g_{t_\theta}}$ is the minimal $G$-invariant character associated to $\theta$. Moreover, any minimal $G$-invariant character  $\mu$ of $N$ is of the form $\widehat{\theta}$ for some $\theta \in \irr N$.

The aim of this survey is to present recent results that have been proved on minimal $G$-invariant characters. Many open questions appear in this framework. In Section \ref{G-conjugacy classes of a normal subgroup}, we recall some remarkable results concerning the structural effect of the $G$-conjugacy classes of $N$. In a dual way, $G$ acts by conjugation on the set (up to isomorphism) of irreducible $\mathbb{K}[N]$-modules for the group algebra $\mathbb{K}[N]$ on a field $\mathbb{K}$, as defined by the well-known Clifford's theorem. This action will allow us to introduce the concept of \textit{minimal $G$-invariant $\mathbb{K}[N]$-modules}. In particular, if $\mathbb{K}=\mathbb{C}$ is the field of complex numbers, this action can be translated into a character-theoretical version. One can then consider the natural action of $G$ on the set of irreducible complex characters of $N$, such that the sum of the elements in each orbit determine the corresponding minimal $G$-invariant character of $N$.

In Section \ref{sección algebras}, we determine the structure of the algebra $\zent{\mathbb{K}[G]} \cap \mathbb{K}[N]$, and, in particular we show the fact that the number of minimal $G$-invariant $\mathbb{K}[N]$-modules (up to isomorphism) coincides with the number of $G$-classes in $N$. Consequently, if $\mathbb{K}=\mathbb{C}$, we obtain that the number of minimal $G$-invariant characters of $N$ is exactly the number of $G$-conjugacy classes of $N$. This fact allows for the definition of the so-called \textit{$G$-invariant table of $N$} (a square and non-singular matrix), as well as the \textit{$G$-character tables of $N$}, which are obtained from the character table of $G$ and contain information on some structural properties of $N$. As we will see, the existence of the $G$-invariant table of $N$ allows to obtain information from the real $G$-conjugacy classes of $N$ (which are not neccesarily constituted from real $N$-conjugacy classes). These results are gathered in Section \ref{G-character tables of a normal subgroup}.

On the other hand, in Section \ref{vanishing elements of minimal G-invariant characters}, we exhibit some results on vanishing and non-vanishing $G$-conjugacy classes of $N$ (that is, $G$-conjugacy classes that are either vanishing or non-vanishing with respect to minimal $G$-invariant characters of $N$). It will be shown that the sizes of these classes affect to some extent the structure of $N$. Moreover, an extension of the well-known Burnside's theorem for characters is presented. Finally, Section \ref{minimal G-invariant character degrees} compiles some known results on the influence of minimal $G$-invariant character degrees of $N$ on its structure. An example of this is an extension of Thompson's theorem, as well as other results regarding the hypercentrality of $N$ as a subgroup of $G$. More precisely, we have that if $N$ possesses only two minimal $G$-invariant character degrees, then $N$ is solvable. We want to remark that, although for $G=N$ this is an elementary result, the proof in this context makes use of the Classification of Finite Simple Groups. In this case (that is, when $N$ has two minimal $G$-invariant character degrees), it is conjectured that the commutator subgroup $[N,G]$ is abelian. We show certain progress that has been made on this problem, which currently remains open. We refer to the original papers for all proofs.

\bigskip

Throughout this paper, the conventions used for group and character-theoretical concepts are those of \cite{Is}. If $G$ is a group and $x\in G$, then $x^G$ is the conjugacy class of $x$ in $G$. The set $\{|x^G| , \ x \in G \}$ is denoted by $\cs G$. Following the notation of \cite{BFM18}, the set $\{|x^G|, \ x \in N\}$ of $G$-conjugacy classes of $N$ is denoted by $\csg G N$.

Moreover, we denote the set of irreducible characters of $G$ by $\irr G$, the set of linear characters of $G$ by $\Lin G$, and the set of irreducible character degrees of $G$ by $\cd G$.  If $N$ is a normal subgroup of $G$ and $\theta$ is an irreducible character of $N$, then we say that $\chi \in \irr G$ \textit{lies over} $\theta$ if $[\chi_N,\theta] \neq 0$, and we denote by
$$\irr {G \mid \theta} = \{\chi \in \irr G \mid [\chi_N , \theta] \neq 0\}$$
the set of irreducible characters of $G$ lying over $\theta$.

\bigskip

\section{$G$-conjugacy classes of a normal subgroup} \label{G-conjugacy classes of a normal subgroup}

Let $G$ be a finite group. If $N$ is a normal subgroup of $G$, then $N$ is a disjoint union of $G$-conjugacy classes. We write
$$N=x_1^G \ \dot{\cup} \cdots \dot{\cup} \ x_r^G,$$
where $x_1, \hdots, x_r \in N$ are representatives of these classes, and where each $x_i^G$, for $1 \leq i \leq r$, can be partitioned into $N$-conjugacy classes,
$$x_i^G=(x_i^{g_1})^N \ \dot{\cup} \cdots \dot{\cup} \ (x_i^{g_{s_i}})^N,$$
with $1=g_1,\hdots,g_{s_i} \in G$. The literature on $G$-conjugacy classes and their influence on the properties of $N$ has predominantly answered the questions of which structure a normal subgroup can have when it contains few $G$-conjugacy classes, when the $G$-conjugacy class sizes are known, or from the (prime, or common divisor) graphs associated to the $G$-conjugacy class sizes. Regarding class sizes, although the elements of $\csg G N$ are multiples of integers in $\cs N$, they may be divisible by primes not dividing $N$, and they do not determine in general the elements of $\cs N$. Deducing information on the conjugacy classes of $N$, unless particular circumstances like a coprime action, usually requires additional, arithmetical hypotheses on $\csg G N$. In a similar way that the set $\cs N$ exerts an influence on $N$, the $\csg G N$ can also affect it. For instance, in \cite{ABF}, it is proved that if all noncentral $G$-conjugacy classes of $N$ have the same size (that is, when $|\csg G N|=2$), then $N$ is nilpotent, and in \cite{ABF14} it is proved that $N$ is solvable when $|\csg G N|=3$. We refer to \cite{BFM18} for an extensive survey on this topic.

In \cite{Be}, it is proved that the $G$-conjugacy class sizes of $N$ and their multiplicities determine whether $N$ is hypercentral in $G$. We recall that a normal subgroup $N$ of $G$ is hypercentral if there exists some integer $r \geq 1$ such that $\textrm{\bf Z}_G^r(N)=N$, where the upper central $G$-series of $N$ is
$$1=\textrm{\bf Z}_G^0(N) \unlhd \zent G \cap N = \textrm{\bf Z}_G^1(N) \unlhd \textrm{\bf Z}_G^2(N) \unlhd \cdots \, ,$$
defined by
$$\textrm{\bf Z}_G^{i+1}(N)/\textrm{\bf Z}_G^i(N)=\zent{G/\textrm{\bf Z}_G^i(N)} \cap N/\textrm{\bf Z}_G^i(N).$$

If $N$ is hypercentral in $G$, then the smallest $r \geq 1$ such that $\textrm{\bf Z}_G^r(N)=N$ is equal to the smallest $r \geq 1$ such that $\Gamma_G^{r+1}(N)=1$, where the lower central $G$-series of $N$ is defined by
$$\Gamma_G^1(N)=N\unrhd \Gamma_G^2(N) = [\Gamma_G^1(N),G] \unrhd \Gamma_G^3(N) = [\Gamma_G^2(N), G] \unrhd \cdots $$
We call the integer $r$ the \textit{hypercentral $G$-length of $N$} (see \cite{AFJP}).

\section{The algebra $\zent{\mathbb{K}[G]}\cap \mathbb{K}[N]$} \label{sección algebras}

Let $\mathbb{K}$ be a field and, for any finite group $G$, denote the group algebra over $\mathbb{K}$ by $\mathbb{K}[G]$. If $N$ is a normal subgroup of $G$, then one may pose the question of which relation exists between the structure of $\mathbb{K}[G]$ and that of $\mathbb{K}[N]$. Particularly, we summarize some results from  \cite{FePRSo} analyzing the structure of $\zent{\mathbb{K}[G]} \cap \mathbb{K}[N]$, which is determined by Clifford's Theorem.

If $V$ is a $\mathbb{K}[G]$-module and $M$ is an irreducible $\mathbb{K}[G]$-module, the $M$-homogeneous component of $V$ is 
$$H_V(M)= \sum \{X \mid X \text{ is a $\mathbb{K}[G]$-submodule of }V,\ V\cong M\},$$
the sum of all those $\mathbb{K}[G]$-submodules of $V$ which are isomorphic to $M$. If $V=V_1\oplus \cdots \oplus V_r$ is completely reducible, then $H_V(M)=\bigoplus\{V_i \mid V_i\cong M\}$. In particular, every completely reducible $\mathbb{K}[G]$-module is a direct sum of its homogeneous components.

Now, let $W$ be a $\mathbb{K}[N]$-module. Let $\overline{W}$ be a copy of the $\mathbb{K}$-vector space $W$, and let $\overline{w}$ denote the image of each $w\in W$ under some $\mathbb{K}$-linear isomorphism from $W$ to $\overline{W}$. For a fixed element $g\in G$, define an action of $N$ on $\overline{W}$ by $$\overline{w}\,n=\overline{w\,(gng^{-1})}$$ for all $w\in W$ and $n\in N$. Then $\overline{W}$ becomes a $\mathbb{K}[N]$-module under this action, called the conjugate module of $W$ by $g$ and denoted by $W^g$ (see Definition (7.2) of Chapter B in \cite{DH}). It easy to see that if $W$ is a $\mathbb{K}[N]$-submodule of a $\mathbb{K}[G]$-module $V$, then $Wg=\{wg\mid w\in W\}\subseteq V$ is a $\mathbb{K}[N]$-submodule of $V$ isomorphic to $W^g$. 

The inertia subgroup of $W$ in $G$ is the subgroup of $G$ defined as $$I_G(W)=\{g\in G\mid W^g\cong W \text{ as } \mathbb{K}[N]\text{-modules}\}.$$ The module $W$ is said to be $G$-invariant if $W^g$ is isomorphic to $W$ for all $g\in G$, \textit{i.e.} if $I_G(W)=G$. If $W$ is an irreducible $\mathbb{K}[N]$-module and $g\in G$, then $W^g$ is also irreducible, and so conjugation defines an action of the group $G$ on the set of irreducible $\mathbb{K}[N]$-modules. The orbit of $W$, up to isomorphism of $\mathbb{K}[N]$-modules, is the set $\{W^{g_i}\mid i=1,\dots,t\}$, where the set $\{g_i\in G\mid i=1,\dots, t\}$ is a right transversal in $G$ of $I_G(W)$. The following result is Lemma (6.4) in \cite{Is}.

\begin{lemma} Let $V$ be a $\mathbb{K}[G]$-module and let $N$ be a normal subgroup of $G$. Suppose $W \subseteq V$ is a $\mathbb{K}[N]$-submodule.

\begin{itemize}
\item[(i)] If $g \in G$, then $Wg$ is a $\mathbb{K}[N]$-module conjugate to $W$.

\item[(ii)] If $M$ is a $\mathbb{K}[N]$-module conjugate to $W$, then $M \cong Wg$ for some $g \in G$.

\item[(iii)] If $U \subseteq V$ is a $\mathbb{K}[N]$-submodule isomorphic to $W$, then $Ug \cong Wg$ as $\mathbb{K}[N]$-modules.
\end{itemize}
\end{lemma}

Let $V$ be a $\mathbb{K}[G]$-module and let $W$ be an irreducible $\mathbb{K}[N]$-module. Denote by $n_W(V)$ the number of $\mathbb{K}[G]$-submodules of $V$ isomorphic to $W$. We now state the Clifford's Theorem (see Theorem 6.5 in \cite{Is}).

\begin{theorem}[Clifford] Let $N$ be a normal subgroup of $G$ and let $V$ be an irreducible $\mathbb{K}[G]$-module. Let $W$ be any irreducible $\mathbb{K}[N]$-submodule of $V$. Then

\begin{itemize}
    \item[(i)] $V=W_1 \oplus \cdots \oplus W_r$, where the $W_i$ are irreducible $\mathbb{K}[N]$-submodules of $V$;

    \item[(ii)] each $W_i$ is of the form $Wg_i$ for some $g_i \in G$ and thus is conjugate to $W$;

    \item[(iii)] we have that $n_W(V)=n_M(V)$ for every $\mathbb{K}[N]$-module $M$ conjugate to $W$.
\end{itemize}
\end{theorem}

The next result is easily checked.

\begin{lemma}\label{conjugate} With the notation above, if $W$ is an irreducible $\mathbb{K}[N]$-module, the following assertions hold.
\begin{enumerate}\item[(i)] The module $\widehat{W}=\bigoplus_{i=1}^tW^{g_i}$ is a minimal $G$-invariant $\mathbb{K}[N]$-module.
\item[(ii)] Every minimal $G$-invariant $\mathbb{K}[N]$-module can be constructed from one orbit of irreducible modules in this way, and this defines a one-to-one correspondence between the set of orbits in the action by conjugation of $G$ on the set of irreducible $\mathbb{K}[N]$-modules and the set of minimal $G$-invariant $\mathbb{K}[N]$-modules (up to isomorphism).
\item[(iii)] Every completely reducible $G$-invariant $\mathbb{K}[N]$-module is a direct sum of  minimal $G$-invariant $\mathbb{K}[N]$-modules.
\end{enumerate}
\end{lemma}

If $e$ is a positive integer and $V$ is any $\mathbb{K}[G]$-module,  we denote by $e\,V=V\oplus\overset{e}{ \cdots}\oplus V$ the direct sum of $e$ copies of $V$.
\begin{lemma}\label{Irr(G)} Let $N$ be a normal subgroup of a group $G$. For irreducible $\mathbb{K}[G]$-modules $V$ and $U$ the following statements are equivalent.
\begin{itemize}\item[(i)] The restricted $\mathbb{K}[N]$-modules $V_N$ and $U_N$, i.e. the modules $V$ and $U$ viewed as $\mathbb{K}[N]$-modules, have a common irreducible $\mathbb{K}[N]$-submodule (up to isomorphism).
\item[(ii)] There exist an irreducible $\mathbb{K}[N]$-submodule $W$ and positive integers $e_V$ and $e_U$ such that $V_N\cong e_V\, \widehat{W}$ and $U_N\cong e_U\, \widehat{W}$.
\item[(iii)] There are positive integers $e_V$ and $e_U$ such that $e_U\, V_N\cong e_V\, U_N$.
\end{itemize}
\end{lemma}

Let $U$ and $V$ be two irreducible $\mathbb{K}[G]$-modules. We say that $U$ and $V$ are \textit{equivalent with respect to $N$} if they satisfy any of the conditions in Lemma \ref{Irr(G)}. This is an equivalence relation where the equivalence class of each $\mathbb{K}[G]$-module is associated to $\widehat{W}$, a minimal $G$-invariant $\mathbb{K}[N]$-module, for some irreducible $\mathbb{K}[N]$-module $W$. Thus, a one-to-one correspondence is defined between the set of equivalence classes of irreducible $\mathbb{K}[G]$-modules and the set of minimal $G$-invariant $\mathbb{K}[N]$-modules (up to isomorphism).

It is easy to prove the following lemma for group algebras in our context. We denote by $\zent{\mathbb{K}[G]}$ the center of $\mathbb{K}[G]$.

\begin{lemma}
\label{lemma_classes}
Let $N$ be a normal subgroup of a group $G$, and $\mathbb{K}$ be a field. Let $\{K_1,\ldots, K_l\}$ be the set of $G$-conjugacy classes of $N$, and let
$\widehat{K_i}$ be the formal sum of $K_i$ in $\mathbb{K}[G]$, for each $i=1,\dots, l$. Then $\{\widehat{K_1},\ldots,\widehat{K_l}\}$ forms a basis of the $\mathbb{K}$-algebra $\zent{\mathbb{K}[G]}\cap \mathbb{K}[N]$.
\end{lemma}

 Let $\mathbb {K}$ denote a splitting field for a group $G$, whose characteristic does not divide the order of $G$, which implies by Maschke's theorem that the group algebra $\mathbb {K}[G]$ is semisimple. 

Since the number of non-isomorphic irreducible $\mathbb{K}[G]$-modules is finite, let the positive integer $k$ denote the number of equivalence classes in the equivalence relation with respect to a normal subgroup $N$, on the set of irreducible $\mathbb {K}[G]$-modules. For each $i=1,\dots, k$, let $\{V_{i1}, \ldots, V_{i{s_i}}\}$, for some positive integer $s_i$, denote a system of pairwise non-isomorphic irreducible $\mathbb{K}[G]$-modules in the same equivalence class, and set
$$T_i=H_{\mathbb{K}[G]}(V_{i1})\oplus \cdots \oplus H_{\mathbb{K}[G]}(V_{i{s_i}}),\ \text{and } L_i=T_i\cap \mathbb{K}[N],$$
which are bilateral ideals of $\mathbb{K}[G]$  and $\mathbb{K}[N]$, respectively.

\begin{theorem} \label{th 3.6}
Let $N$ be a normal subgroup of a group $G$. The following statements hold:
\begin{enumerate}
\item[(i)] $\mathbb{K}[N] = L_1 \oplus \cdots \oplus L_k$;
\item[(ii)] $\zent{\mathbb{K}[G]}\cap \mathbb{K}[N] = (L_1\cap \zent{T_1}) \oplus \cdots \oplus (L_k\cap \zent{T_k})$;
\item[(iii)] $dim_{\mathbb{K}}(L_i\cap \zent{T_i})=1$, for every $i=1,\dots,k$;
\item[(iv)] $dim_{\mathbb{K}}(\zent{\mathbb{K}[G]}\cap \mathbb{K}[N])=k$, which is equal to the number of non-isomorphic minimal $G$-invariant $\mathbb{K}[N]$-modules.
\end{enumerate}
\end{theorem}

As a consequence of Lemma \ref{lemma_classes} and Theorem \ref{th 3.6}, we obtain the following.

\begin{corollary} The number of orbits in the action by conjugation of $G$ on the set of irreducible $\mathbb{K}[N]$-modules (that is,  the number of non-isomorphic minimal $G$-invariant $\mathbb{K}[N]$-modules) is equal to the number of $G$-conjugacy classes of $N$.

\end{corollary}

\section{$G$-character tables of a normal subgroup} \label{G-character tables of a normal subgroup}

In the rest of the paper we focus on characters of groups, which are always considered over the field $\mathbb{K}=\mathbb{C}$ of complex numbers. Then, all results of Section \ref{sección algebras} can be applied.

The character-theoretical version of Clifford's Theorem (see Theorem (6.2) of \cite{Is}), as well as Clifford's correspondence for characters (see Theorem (6.11) of \cite{Is}) are essential in the study of characters of normal subgroups of finite groups. Let $\theta \in \irr N$ and suppose $\chi \in \irr {G \mid \theta}$. Clifford's theorem states that 
$$\chi_N=e \widehat{\theta},$$
where the integer $e$ is called the \textit{ramification number}, and
$$\widehat{\theta}=\theta^{g_1}+\cdots+\theta^{g_{t_\theta}}$$
is the minimal $G$-invariant character associated to $\theta$.

We have that $t_\theta$ is the index in $G$ of $I_G(\theta)$, the inertia subgroup of $\theta$ in $G$. We denote the set of minimal $G$-invariant characters of $N$ by
$$ \Min G N=\{\widehat{\theta}\mid \theta \in \irr{N}\}.$$

The following theorem, which follows from translating Section \ref{sección algebras} to characters, can also be deduced from Corollary (6.33) of \cite{Is}.

\begin{theorem}\label{brauer} Let $N$ be a normal subgroup of the group $G$. Then the number of orbits in the action by conjugation of $G$ on $\irr N$ and on the set of conjugacy classes of $N$ are equal. That is, the number of $G$-conjugacy classes of $N$ coincides with the number of elements in $\Min G N$.
\end{theorem}

Two characters $\chi, \psi \in \irr{ G }$ can be said to be \textit{equivalent with respect to $\theta$} if there exists $\theta \in \irr N$ such that $[\chi_N, \theta] \neq 0$ and $[\psi_N, \theta] \neq 0$. This is indeed an equivalence relation on $\irr G$
where each equivalence class is the set $\irr{G |\theta}$ for a $\theta \in \irr N$. More precisely, on account of Clifford's theorem, we have the following lemma.

\begin{lemma}\label{Irr(G)-characters} Let $N$ be a normal subgroup of a group $G$. For $\chi,\psi \in \irr{G}$ the following statements are equivalent:
\begin{itemize}\item[(i)] $[\chi_N,\psi_N]\neq 0$.
\item[(ii)] There exist $\theta\in \irr{N}$ and positive integers $e_\chi$ and $e_\psi$ such that $\chi_N= e_\chi\, \widehat{\theta}$ and $\psi_N= e_\psi\, \widehat{\theta}$.
\item[(iii)] There is a rational number $c$  such that $\chi_N= c\, \psi_N$.
\end{itemize}
\end{lemma}

This equivalence relation defines a one-to-one correspondence between the set of equivalence classes of irreducible characters of $G$ and the set of minimal $G$-invariant characters of $N$.

It is now possible to introduce the following concept: a \textit{$G$-character table of $N$} is any matrix $\textsf{X}=(x_{ij}) \in \textup{M}_k(\mathbb{C})$ with entries $$x_{ij}=\chi_i(n_j),\quad 1\leq i,j \leq k, $$ where  $\{n_1^G, \ldots, n_k^G\}$  is the set of $G$-conjugacy classes of $N$ and $\Delta = \{\chi_1=1_G, \chi_2, \dots, \chi_k\}$ denotes a representative system of the equivalence classes in the equivalence relation with respect to $N$, defined on $\irr{G}$. By Theorem \ref{brauer}, we have that $G$-character tables $\textsf{X}$ of $N$ are square.

\medskip

We now fix some notation:
\begin{itemize}\item  $\{n_1^G, \ldots, n_k^G\}$  the set of $G$-conjugacy classes of $N$;
 \item $\textsf{D}=(d_{ij})$ a diagonal matrix with entries $d_{ij}=\delta_{ij}|n_i^G|$, where $\delta_{ij}$ is the Kronecker delta function, for $i,j\in \{1, \ldots, k\}$;
\item $\Delta = \{\chi_1=1_G, \chi_2, \dots, \chi_k\}$ a representative system of the equivalence classes in the equivalence relation with respect to $N$, defined on $\irr{G}$;
\item $\Omega= \{\theta_i \in \irr{N} \mid 1 \leq i \leq k\}$ such that $\chi_i \in \irr{G\mid \theta_i}$, for each $i \in \{1, \dots, k\}$;
\item $t_i=t_{\theta_i}=|G:I_{G}(\theta_i)|$, $e_{i} = [(\chi_i)_N, \theta_i]\neq 0$, for $1\le i\le k$;
\item $\textsf{X}$ the $G$-character table of $N$ constructed from $\Delta$, and $\overline{\textsf{X}}^t$ its transposed conjugate matrix; i.e if $\textsf{X}=(x_{ij})$, then $\overline{\textsf{X}}^t=(y_{ij})$ being $y_{ij}=\overline{x_{ji}}$ the complex conjugate of $x_{ji}$, for $1\le i,j\le k$;
\item  $\Lambda_{\textsf{X}}=\text{diag}(\lambda_1,\dots,\lambda_k)$ the  diagonal matrix with entries $\lambda_i=|N|t_{i} e_i^2$, for $ 1\le i\le k.$
\end{itemize}

\medskip

The results below are proved in \cite{FePRSo}.

\begin{theorem}\label{square}
We have that  \textup{$\Lambda_{\textsf{X}}=\textsf{X}\textsf{D}\overline{\textsf{X}}^t.$} In particular, \textup{$\textsf{X}$} is non-singular.
\end{theorem}

\begin{corollary}\label{restriction} Let $\chi\in \irr{G}$ and $i\in \{1,\dots,k\}$ such that $\chi=\chi_i\in \Delta$. Then $\chi_N\in \irr{N}$ if and only if  $\lambda _i=|N|$.
\end{corollary}

From  Theorem~\ref{square}, one can deduce that some arithmetical relations between the aforementioned integers $t_i$, $e_{i}$ and $\theta_i(1)$ can be read off the matrix $\Lambda_{\textsf{X}}$, for the  $G$-character table $\textsf{X}$ of $N$.

\begin{corollary}\label{relations}
 The next integer relations hold, where the corresponding right sides can be computed from the  $G$-character table $\textsf{X}$ of $N$:
$$e_{i}^2t_i=\dfrac{\lambda_i}{|N|}  $$
$$t_i\theta_i(1)^2=\dfrac{|N|\chi_i(1)^2}{\lambda_i} $$
$$\dfrac{\theta_i(1)}{e_{i}}=\dfrac{|N|\chi_i(1)}{\lambda_i}$$
$1\le i\le k$.
\end{corollary}

Choosing a different representative system for the equivalence classes in the equivalence relation on $\irr G$ with respect to $N$ may yield a different $G$-character table. Nevertheless, the corresponding rows in two $G$-character tables of $N$ are always proportional, by a rational factor given by the ramification numbers. In order to obtain a table that is unique, we can define the \textit{$G$-invariant table of $N$} as the matrix $\widehat{\textsf{X}}= (y_{ij}) \in \textup{M}_k(\mathbb{C}),$
with entries $$y_{ij}=\widehat{\theta}_i(n_j),\quad 1\leq i,j \leq k. $$
Note, however, that this table cannot be read from the character table of $G$, whereas the $G$-character tables of $N$ can be. As the rows of the $G$-invariant table of $N$ only differ from the corresponding ones in a $G$-character table of $N$ by a rational multiplicative factor, it is also a non-singular matrix, by Theorem \ref{square}.

From the $G$-invariant table of $N$, we obtain the following three results about real $G$-conjugacy classes.

\begin{lemma}[\cite{FePRSo}, Corollary 5.4] Let $N$ be a normal subgroup of a group $G$. Then the number of real $G$-conjugacy classes in $N$ is equal to the number of real-valued minimal $G$-invariant characters of $N$.
\end{lemma}

\begin{lemma}[\cite{FePRSo}, Lemma 5.6] Let $N$ be a normal subgroup of a group $G$. If there is a unique real $G$-conjugacy class in $N$, then $N$ is of odd order, but the converse is not true in general.
\end{lemma}

\begin{lemma}[\cite{FePRSo}, Corollary 5.7] Let $N$ be a normal subgroup of a group $G$. Assume that there are exactly two real $G$-conjugacy classes in $N$. Then $N$ has a normal Sylow $2$-subgroup.
\end{lemma}

One might wonder, from the previous lemmas, whether a relation exists between the rational $G$-classes of $N$ and its minimal $G$-invariant characters with rational values.

\section{Vanishing elements of minimal $G$-invariant characters} \label{vanishing elements of minimal G-invariant characters}

It was seen previously that counting the real-valued minimal $G$-invariant characters in the $G$-character tables of $N$ can provide some information on the normal subgroup. We now focus on the zero values in $G$-character tables. We say that $x \in N$ is \textit{vanishing} in $G$ if there exists some $\widehat{\theta} \in \Min G N$ such that $\widehat{\theta}(x)=0$. That is, if there is $\theta \in \irr N$ that vanishes on $x$. On the contrary, we say that $x$ is \textit{non-vanishing} if no $\widehat{\theta} \in \Min G N$ vanishes on $x$. These concepts can be applied without ambiguity to the conjugacy class of $x$.

In this section, we show some known results which evidence that vanishing  or non-vanishing  $G$-conjugacy classes of minimal $G$-invariant characters also exert a certain influence on the structure of a normal subgroup. In particular, the following extension of Burnside's theorem can be useful in order to obtain vanishing elements of $N$ in $G$ (see Section 2 in \cite{FGS26}).

\begin{proposition} \label{burnside}
    Let $N$ be a normal subgroup of $G$ and let $\widehat{\theta} \in \Min G N$. Let $K=g^G$ be the $G$-conjugacy class of an element $g \in N$. If $(\widehat{\theta}(1), |K|) = 1$, then either $g \in \zent{\widehat{\theta}}$ or $\widehat{\theta}(g) = 0$.
\end{proposition}

Note that Proposition \ref{burnside} cannot be obtained as a consequence of Burnside's theorem for ordinary characters. As a matter of fact, it cannot be applied to a $\chi \in \irr{G \mid \theta}$, since the hypothesis $(\widehat{\theta}(1),|K|)=1$ does not imply that $(\chi(1),|K|)=1$. On the other hand, if $(\widehat{\theta}(1),|K|)=1$, then $(\theta(1),|g^N|)=1$, for $g^N$ the conjugacy class of $g$ in $N$, then Burnside's theorem applied to $\theta$ guarantees that $g \in \zent \theta$ or $\theta(g)=0$, but this does not imply that $g \in \zent{\widehat{\theta}}$ or $\widehat{\theta}(g)=0$.

It is easy to deduce from Proposition \ref{burnside} the following result about $G$-conjugacy class sizes of $N$ which are of prime power size.

\begin{corollary}
Let $N$ be a minimal normal subgroup of $G$. Suppose that $K=x^G$ is the $G$-conjugacy class of a non-central element $g \in N$ such that $|K|$ is a power of a prime $p$. Then $p$ divides $|G:N|$.
\end{corollary}

The vanishing $G$-conjugacy class sizes also contain information on certain properties of the normal subgroup, as seen below.

\begin{theorem}[\cite{FGS20}, Theorem C]
Let $N$ be a normal subgroup of a group $G$, and let $\pi$ be any set of prime numbers.
\begin{enumerate}
\item[(i)] Suppose that $|x^G|$ is a $\pi'$-number for every prime-power order $\pi$-element $x \in N$ which is vanishing in $G$. If $N$ is $\pi$-separable, then $N/O_{\pi'}(N)$ has a nilpotent normal Hall $\pi$-subgroup. In particular, the Hall $\pi$-subgroups of $N$ are nilpotent.
    
\item[(ii)] Suppose that $|x^G|$ is a $\pi$-number for every prime-power order $\pi'$-element $x \in N$ which is vanishing in $G$. If $\text{Hall}_{\pi}(N) \neq \emptyset$, then $N$ has a normal Hall $\pi$-subgroup. Additionally, if all $|x^G|$ are also $\pi'$-numbers for the prime-power order $\pi$-elements $x \in N$ that are vanishing in $G$, then the Hall $\pi$-subgroups of $N$ are nilpotent.
\end{enumerate}
\end{theorem}

In the opposite way, it is possible to deduce some properties of $N$ by studying the non-vanishing elements of minimal $G$-invariant characters. Namely, if any $G$-character table of $N$ has no zeros on elements of prime-power order, then $N$ is nilpotent.

\begin{theorem}[\cite{FGS20}, Theorem A] \label{vanishing 1}
Let $N$ be a normal subgroup of a group $G$, and let $p$ be a prime. If 
$\widehat{\theta}(x) \neq 0$,
for every $p$-element $x \in N$ and for all $\theta \in \irr N$, then $N$ has a normal Sylow $p$-subgroup.  In particular, if $\widehat{\theta}(x) \neq 0$ for every prime power-order element $x \in N$ and for all $\theta \in \irr N$, then $N$ is nilpotent.
\end{theorem}

The theorem above generalizes for $N=G$ a previous result of S. Dolfi, E. Pacifici, L. Sanus and P. Spiga (see \cite{DPSS}). We remark that a direct application of their theorem does not prove Theorem \ref{vanishing 1}, as $x$ being a non-vanishing element for $\theta$ does not imply that it is non-vanishing for $\widehat{\theta}$, nor vice versa. 

In fact, the normality of the Sylow $p$-subgroups of $N$ can be characterized through the minimal $G$-invariant characters of $N$, as seen in the theorem below. This, in turn, generalizes a result by G. Malle and G. Navarro for $N=G$ (see Theorem 3.3 of \cite{MN12}).

\begin{theorem}[\cite{FGS20}, Theorem B] Let $N$ be a normal subgroup of a group $G$, and let $P$ be a Sylow $p$-subgroup
of $G$ for some prime $p$. Let $P_0 = P \cap N$ and $\beta \in \irr{P/P_0}$. 
Then the following conditions are pairwise equivalent:
\begin{enumerate}
    \item[(i)] $P_0$ is a normal Sylow $p$-subgroup of $N$;
    \item[(ii)] $\widehat{\theta}(x) \neq 0$ for all irreducible constituents $\theta$ of $(1_{P_0})^N$ and all $x \in P_0$\,;
    \item[(iii)] $\widehat{\theta}(x) \neq 0$ for all $\theta \in {\rm Irr}(N/P_0)$ such that $[\theta, \beta^G_{N}] \neq 0$ and all $x \in P_0$\,.
\end{enumerate}
\end{theorem}

\section{Minimal $G$-invariant character degrees} \label{minimal G-invariant character degrees}

In this last section, we focus on the degrees of minimal $G$-invariant characters. It is natural to wonder whether they may provide information on the structure of $N$, in the same way that $\cd G$, the set of irreducible character degrees of $G$, can do for the structure of $G$. We denote by
$$\Mcd G N = \{\widehat{\theta}(1)\mid \theta \in \irr{N}\}$$
the set of \textit{minimal $G$-invariant character degrees}, and by
$$\text{Lin}_G(N)=\{\widehat{\theta} \in \Min G N \mid \widehat{\theta}=\theta \text{ and } \theta(1)=1\},$$
the set of \textit{linear $G$-invariant characters} of $N$.

The lemma below determines the characters in $\text{Lin}_G(N)$. Moreover, it shows $\Mcd G N = \{1\}$ if and only if $N$ is central in $G$. 

\begin{lemma}[\cite{AFJP}, Proposition 3.2] Let $N$ be a normal subgroup of a group $G$ and let $\theta \in \Irr(N)$. Then $$[N,G]=\bigcap_{\theta \in \lin G N}\ker \theta.$$
In particular, $|\,N:[N,G]\,|=|\,\lin G N\,|$.
\end{lemma}

Observe that if none of the minimal $G$-invariant character degrees of $N$ is divisible by a certain prime $p$, then trivially $N$ has an abelian normal Sylow $p$-subgroup by Itô-Michler's theorem. Conversely, the following theorem extends Thompson's theorem to minimal $G$-invariant character degrees.

\begin{theorem}[\cite{AFJP}, Theorem 3.10]
Let $N$ be a normal subgroup of the group $G$ and let $p$ be a prime integer. Suppose $p$ divides $\widehat{\theta}(1)$ for every $\widehat{\theta}\in \Min G N \setminus \lin G N$. Then $N$ has a normal $p$-complement.
\end{theorem}

Specifically, when all the elements in $\Mcd G N$ are  $p$-powers for a prime $p$, we have the following lemma. Note that the condition of $N$ being a $p$-group is necessary (see Example 2 in \cite{FGS26}).

\begin{lemma}[\cite{FGS26}, Proposition 3.2 ; \cite{AFJP}, Corollary 3.9]\label{hypercentral}
Let $N$ be a normal subgroup of the group $G$. If $N$ is a $p$-group for some prime $p$ and $\widehat{\theta}(1)$ is a $p$-power number for all $\widehat{\theta} \in \Min G N$, then $N$ is a hypercentral subgroup of $G$. In particular, $[N, G]\subseteq {\bf \Phi}(G)$, where $\Phi(G)$ denotes the Frattini subgroup of $G$.
\end{lemma}

In the ordinary case, we have that a group is solvable if it possesses only two irreducible character degrees, and this can be proved through elementary methods. In our context, we have that a normal subgroup $N$ with two minimal $G$-invariant character degrees is also solvable, but the proof of this fact requires the CFSG.

\begin{theorem}[\cite{FGS26}, Theorem A]\label{2 minimales}
Let $N$ be a normal subgroup of the group $G$. If $\Mcd G N = \{1, f\}$, for some positive integer $f$, then $N$ is solvable.
\end{theorem}

When $N$ has three minimal $G$-invariant character degrees, an open question is whether $N$ is still solvable, as it happens for irreducible character degrees when $N=G$.

On the other hand, Corollary (12.6) of \cite{Is} states that the derived subgroup of a group with only two irreducible character degrees is abelian. In the same direction, it is conjectured (Conjecture B of \cite{FGS26}) that, if $|\Mcd G N|=2$, then $[N,G]$ is also abelian. Although proving this is still an open problem, we have that if $|\Mcd G N|=2$ and $[N,G]$ is not abelian, then the degree of the characters in $\Min G N \setminus \text{Lin}_G(N)$ is a power of $p$. In that case, when $N$ is a $p$-group, we have that $N$ is hypercentral in $G$ from Lemma \ref{hypercentral}.

\begin{theorem}[\cite{FGS26}, Theorem C]
Let $N$ be a normal subgroup of the group $G$. If $\Mcd{G}{N}=\{1,f\}$, then either $[N,G]$ is abelian or $f$ is a power of $p$, for some prime $p$.
\end{theorem}

One result that may be useful to solve the mentioned conjecture is the following generalization of Taketa's theorem. First, we need to introduce some definitions. The characters $\theta \widehat{\theta}$ of $N$ are called the \textit{leader $G$-characters} of $N$. We denote by
$$\Lcd G N = \{(\theta\widehat{\theta})(1)\mid \theta \in \irr{N}\}$$
the set of \textit{leader $G$-character degrees} of $N$. A normal subgroup $N$ is said to be a \textit{$G$-invariant nMI-subgroups} of $G$ if for every $\theta \in \irr{N}$, there exist $H_{\theta} \unlhd G$, with $H_{\theta} \subseteq N$, and $\lambda_{\theta} \in \text{Lin}_G(H_{\theta})$ such that $\irr{G \mid \theta}=\irr{G \mid \lambda_{\theta} }$ (that is, every irreducible character of $N$ is fully ramified with respect to a linear $G$-invariant character of a $G$-invariant normal subgroup of $N$).

 \begin{theorem}[\cite{FGS26}, Theorem 2.9] \label{Lcd}
 Let $N$ be a $G$-invariant $nMI$-subgroup of a group $G$ and let $1=f_1 < f_2 < \cdots < f_s$  the distinct elements in $\textnormal{Lcd}_G(N)$. Then
$$\Gamma^{i+1}_G(N) \subseteq \ker {\widehat{\theta_{i}}}$$
with $(\theta_i\widehat{\theta_i})(1) = f_i$ and $\theta_i \in \irr{N}$, for $i \in \{1, \ldots, s\}$. In particular, $N$ is a hypercentral subgroup of $G$ and $l_G(N) \leq s=|\Lcd{G}{N}|$.
\end{theorem}

We note that Theorem \ref{Lcd} does not hold if we replace the set $\Lcd G N$ by $\Mcd G N$ (see Example 3.23 in \cite{AFJP}). By using Theorem \ref{Lcd}, we obtain that $[N,G]$ is abelian for $f=p$ in Theorem \ref{2 minimales}. This suggests that the conjecture of \cite{FGS26} could be true.

\begin{theorem}[\cite{FGS26}, Theorem D]
Let $N$ be a normal subgroup of the group $G$. If $\Mcd G N = \{1, p\}$, for some prime number $p$, then $[N,G]$ is abelian.
\end{theorem}

We finish this section by observing that, in the same way that there exist graphs for the $G$-conjugacy class sizes of $N$, a new common-divisor graph (or, from a dual perspective, a new prime graph) could be defined on the set $\Mcd G N \setminus \{1\}$. From the results presented in this section, we obtain information about the structure of a normal subgroup $N$ of a group $G$ when both new graphs reduce to a single vertex. However, this is a reformulation on minimal $G$-invariant characters that remains open.

\bigskip

\noindent \footnotesize{\textsc{Mar\'{\i}a Jos\'e Felipe, Institut Universitari de Matem\`atica Pura i Aplicada, \newline
Universitat Polit\`ecnica de Val\`encia, 
Val\`encia, Spain.} \newline
\texttt{mfelipe@mat.upv.es}}

\bigskip

\noindent \footnotesize{\textsc{Iris Gilabert, Institut Universitari de Matem\`atica Pura i Aplicada, \newline
Universitat Polit\`ecnica de Val\`encia, 
Val\`encia, Spain.} \newline
\texttt{igilman@posgrado.upv.es}}

\bigskip

\noindent \footnotesize{\textsc{Lucia Sanus, Departament de Matem\`atiques, Facultat de
 Matem\`atiques, \newline
Universitat de Val\`encia,
46100 Burjassot, Val\`encia, Spain.} \newline
\texttt{lucia.sanus@uv.es}}


\begin{thebibliography}{99}

\bibitem{ABF14} Z. Akhlaghi, A. Beltr\'an, M. J. Felipe, Normal sections, class sizes and solvability of finite groups, J. Algebra, 399, 220-231 (2014).

\bibitem{AFJP} Z. Akhlaghi, M. J. Felipe, M. K. Jean-Philippe, Some properties of normal subgroups determined from character tables, Bull. Malays. Math. Sci. Soc., 47(90) (2024).

\bibitem{ABF} E. Alemany, A. Beltr\'an, M. J. Felipe, Nilpotency of normal subgroups having two G-class sizes, Proc. Amer. Math. Soc., 138, 2663–2669 (2011).

\bibitem{Be} A. Beltr\'an, A criterion for a normal subgroup to be hypercentral based on class sizes, Ricerche Mat., 74, 143–150 (2025).

\bibitem{BFM18} A. Beltr\'an, M. J. Felipe, C. Melchor, Conjugacy classes contained in normal subgroups: an overview, Int. J. Group Theory, 7 (1), 23-26 (2018).

\bibitem{DH} K. Doerk, and T. Hawkes,	\textit{Finite soluble groups}, De Gruyter (1992).

\bibitem{DPSS}  S. Dolfi, E. Pacifici, L. Sanus, P. Spiga, On the orders of zeros of irreducible characters. J. Algebra, 321, 345–352 (2009).

\bibitem{FGS20} M. J. Felipe, N. Grittini, V. Sotomayor, On zeros of irreducible characters lying in a normal subgroups, Ann. Mat. Pura Appl., 199, 1777-1787 (2020).

\bibitem{FGS26} M. J. Felipe, I. Gilabert, L. Sanus, On degrees of minimal invariant characters, preprint (2026).

\bibitem{FePRSo}  M. J. Felipe, M. D. Pérez-Ramos, V. Sotomayor, On $G$-character tables for normal subgroups, Quaest. Math., 46(4), 721–743 (2022).

\bibitem{Is} I. M. Isaacs, \textit{Character theory of finite groups}. Corrected reprint of the 1976 original. AMS Chelsea Publishing, Providence, RI, 2006.

\bibitem{MN12} G. Malle, G. Navarro, Characterizing normal Sylow $p$-subgroups by character degrees. J. Algebra, 370, 402–406 (2012).
\end{thebibliography}
\end{document}